       \documentclass[11pt]{article} %,openbib
\usepackage{amsmath, amsfonts, mathrsfs, textcomp,amssymb,natbib}
\usepackage[dvips]{graphicx}
         \newenvironment{reflist}{\begin{list}{}
         {\itemsep=20pt \parsep=3pt
          \topsep=0pt  \parskip=20pt  \listparindent=-.15in
           \leftmargin= 0.15in }
         \item \ \vspace{-.35in} }
         {\end{list}}

\newtheorem{theorem}{Theorem}

\newtheorem{proposition}{Proposition}
\newcommand{\ninf}{\mbox{ as } n\to\infty}

\newcommand{\szer}{\mbox{ as } s\to 0}
\newcommand{\tzer}{\mbox{ as } t\to 0}

\newcommand{\xder}[1]{ X^{(#1)}}

\newcommand{\pcal}{{\cal P }}
\newcommand{\bcal}{{\cal B }}
\newcommand{\ccal}{{\cal C }}

\newcommand{\rcal}{{\cal R }}

\newcommand{\norm}[1]{||\, #1\,||}
\newcommand{\oo}{\mbox{o}}

\newcommand{\bt}{\beta}
\newcommand{\ah}{\alpha}

\newcommand{\Ep}{\mbox{E}}

\newcommand{\eps}{\epsilon}
\newcommand{\dt}{\delta}
\newcommand{\ot}{\bar{t}}
\newcommand{\os}{\bar{s}}

\newcommand{\beq}{ \begin{equation}}
\newcommand{\eeq}{ \end{equation}}
\newcommand{\beqr}{ \begin{eqnarray}}
\newcommand{\eeqr}{ \end{eqnarray}}
\newcommand{\beqrn}{ \begin{eqnarray*}}
\newcommand{\eeqrn}{ \end{eqnarray*}}
\newcommand{\bye}{\end{document}}

  \title{Spline approximation of a random process with singularity}
\author{Konrad Abramowicz,$\quad$ Oleg Seleznjev,\\
 Department of Mathematics and Mathematical Statistics\\
Ume{\aa } University, SE-901 87 Ume\aa , Sweden }
        \date{}
\begin{document}

 \maketitle

\begin{abstract}
Let a continuous random process $X$ defined on $[0,1]$ be $(m+\beta)$-smooth, $0\le m$, $0<\beta\le 1$, in quadratic mean for all $t>0$  and have an isolated
singularity point at $t=0$. In addition, let $X$  be locally like a $m$-fold integrated $\beta$-fractional Brownian motion for all nonsingular points. We consider
approximation of $X$ by piecewise Hermite interpolation splines with $n$ free knots (i.e., a \emph{sampling design},  a \emph{mesh}).
The approximation performance is measured by mean errors (e.g., integrated or maximal quadratic mean errors). We
construct a sequence of sampling designs with asymptotic
approximation rate $n^{-(m+\beta)}$ for the whole interval.

\end{abstract}

\textbf{Keywords}: Approximation, random process, sampling design, Hermite splines

 \baselineskip=3.4 ex

\section{Introduction}
Let a random  process $X(t),\,t\in[0,1]$,  with finite second moment be observed in $n$ points and have different  quadratic mean (q.m.)\
smoothness at $t=0$ (i.e., an isolated \emph{singularity}). We consider  piecewise polynomial approximation of $X$, combining two different
Hermite interpolation splines for the interval adjacent to the singularity point  and for the remaining part (i.e., a \emph{composite
Hermite interpolation spline}). A sequence of sampling designs (i.e., meshes) is constructed  to improve the asymptotic approximation
performance (e.g., integrated or maximal q.m.\ errors) eliminating the effect of the singularity point. The proposed technique can  be
applied also  to random processes with a  finite number of isolated singularity points. In principle, the main idea is well-known for
nonlinear approximation of deterministic functions (see, e.g., de Boor, 1973,  DeVore, 1998, and references therein). We use a finer mesh
where the target function is singular and a coarser mesh where it is smooth. The primary question remains however: how to measure this
smoothness for random functions in order to obtain definitive results.

More precisely, let $l$-th   derivative of $X$ satisfy an $\alpha$-H\"older condition $0\le l, 0<\alpha \le 1$, for $[0,1]$.
It is known that the approximation rate  $n^{-(l+\alpha)}$ is optimal in a certain sense for such class  of processes and Hermite splines
(see, e.g., Buslaev and Seleznjev, 1999, Seleznjev, 2000).
%\citep[see][]{Buslaev and Seleznjev, 1999, Seleznjev, 2000}.
Let, additionally, $X$ have a continuous $m$-th derivative satisfying a local stationarity condition, Berman (1974), with parameter $\beta$, and $0\le m, 0<\beta < 1, l+\alpha\leq m+\beta$,
 for all  points $t\in (0,1]$.
We approximate $X$ by  a corresponding composite Hermite interpolation spline and set up a  sequence  of sampling designs attaining   the  asymptotic
 approximation rate  $n^{-(m+\beta)}$ for the whole interval $[0,1]$.

The nonlinear approximation results   obtained in this  paper  are related to various problems in signal processing (e.g., optimization
of compressing digitized signals, see, e.g., Cohen and D'Ales, 1997, Cohen et al., 2002),  in numerical analysis of random functions
(see, e.g., Benhenni and Cambanis, 1992,   Kon and Plaskota, 2005,  Creutzig and Lifshits, 2006, Creutzig et al., 2007), in simulation studies with
controlled accuracy for functionals on realizations of random
processes (see, e.g., Eplett, 1986, Abramowicz  and Seleznjev, 2008).
 Hermite spline interpolation for continuous and smooth random functions
is studied in Seleznjev (1996, 2000), and Huesler et al.\ (2003). Ritter (2000) contains a very detailed survey  of various random function approximation problems.

The paper is organized as follows.
First we introduce a basic notation.
In Section \ref{se:Res}, we consider composite spline interpolation for certain H\"older type classes of random processes, which behave
locally like $m$-fold integrated fractional Brownian motion for all intervals $[a,b]\subset (0,1]$ and construct  sequences of sampling
designs with asymptotically optimal approximation rate. In the second part of the section, the approximation of more general
H\"older's classes of random functions and the accuracy for the corresponding composite Hermite spline approximation are studied.
Sections   \ref{se:Exp} and \ref{se:Pro} contain the results of numerical experiments and
 the proofs of the statements from Section \ref{se:Res}, respectively.

\subsection{Basic notation}
Let $X=X(t),\,t\in[0,1]$, be defined on a probability space $(\Omega,\mathscr{F},P)$. Assume that for every $t$, the random variable $X(t)$
lies in the normed linear space $L^2(\Omega)=L^2(\Omega,\mathscr{F},P)$ of random variables with finite second moments and identified
equivalent elements with respect to $P$.  We set
 $\norm{Z}:=(\Ep\,Z^2)^{1/2}$ for all $Z \in L^2(\Omega)$. Let  $\ccal^m([0,1])$  denote the space of random processes with continuous
q.m.\
derivatives up to order $m\ge 0$, and $C^{m}([0,1])$ be the corresponding space of non-random functions.
Let $\pcal_k([0,1])$ be  the space of stochastic polynomials of order $k\geq 0$.
We define  the norm for any $X\in\ccal^m([0,1])$ by setting
 $$
  \norm{X}_p :=\left(\int^1_0\norm{X(t)}^p dt \right)^{1/p},
   \quad 1\le p <
 \infty,
 $$
 and   $\norm{X}_\infty :=\max_{[0,1]} \norm{X(t)}$ for $p=\infty$.
 Henceforth, we use the convention, that if $p=\infty$, then $1/p=0$.
 Denote by $f(t) \asymp g(t)$ as $t \to a$ the weak equivalence property, i.e.,  $ C_1~g(t)~\le~f(t)~\le~C_2~g(t)$ for some positive $C_1,C_2, b$  and for all $t\in (a,b]$
 and similarly for $a_n\asymp b_n$ as $n\to \infty$.

Let $X$ be sampled at the distinct design points $T_n   :=(t_0,  t_1,\ldots, t_n) $  (also referred to
  as {\it knots}),
  and  the set of all $(n+1)$-point designs be denoted
  by $D_n:=\{T_n:0=t_0<t_1<\cdots<t_n=1\} $. We
   suppress  the argument $n$
  for the  design points $t_k=t_k(n)$ from   $T_n$, $\; k=0,\ldots,n$,
   when doing so causes no confusion.
Recall that,
  for any $f\in C^m([0,1]), \;m\ge 0$, the piecewise
  Hermite polynomial $H_{k}(t)=H_{k}(f,T_n)(t)$, of degree  $k=2r+1$, $0\le r\le  m$,
  is the unique
  solution of the interpolation problem
  $H_{k}^{(j)}(t_i)=f^{(j)}(t_i)$, where
    $i=0,\ldots,n$, $ j=0,\ldots,r$.
Analogously, we suppose that for  $X\in\ccal^{m}([0,1])$, the process and    its first  $r\leq m$ derivatives can be sampled, and write
$H_{k}(X,T_n)$ with $k=2r+1\leq 2m+1$ to denote a corresponding stochastic Hermite spline. Define ${H}_{q,k}(X,T_n),\,q \leq k$, to be a {\it composite
Hermite spline}
$$
{H}_{q,k}(X,T_n) :=
\left\{
\begin{array}{rl}
H_q(X,T_n)(t), & t\in[0,t_1]\\
H_k(X,T_n)(t), & t\in[t_1,1]
\end{array}
\right. .
$$
To formulate the asymptotic results, we introduce a {\it quasi regular
 sequence} (qRS) of sampling designs $\{T_n=T_n(h)\}$
   generated by a {\it positive}  density function $h(t),\,t\in(0,1]$, via
 \beq \label{defds}
   \int_0^{t_{i}} h(t)\;dt =i/n,\quad i=1,\ldots, n,\quad \int_0^1 h(t)dt=1.
   \eeq
  Henceforth we assume that $h(t)$ is continuous  for $t\in(0,1]$ and if $h(t)$ is  unbounded in $t=0$, then $ h(t)\to +\infty$ as $t\to
  0+$, otherwise $h\in C([0,1])$.
We denote this property of  $\{T_n\}$ by: {\em $\{T_n\}$ is qRS$(h)$}.
If the function $h$ is positive and continuous on the whole interval $[0,1]$,  the corresponding generated sampling designs are
called regular sequences,
RS$(h)$, and were introduced by  Sacks and Ylvisaker (1966)  for  certain  time series models.
In particular, if $h$ is uniform over $[0,1]$  ($h(t)\equiv 1,\,t\in[0,1]$), then the  regular sampling becomes the equidistant
sampling  including the endpoints.
For random process linear approximation problems, see Su and Cambanis (1993), Seleznjev (2000), and references therein.
Define the related distribution functions
 $$
 H(t):=\int_0^{t}h(u)du,\qquad G(t):=H^{-1}(t)=\int_0^{t}g(v)dv,\quad t\in[0,1],
 $$
i.e., $G$ is a quantile function for the distribution $H$. Then by the definition,
 $$t_i=G\left({i}/{n}\right),\qquad i=1,\ldots, n,
 $$
the knots are $i/n$-quantile points of  $H$.
The \emph{quantile density function}
$$
g(s):=G'(s)={1}/{h(G(s))}
$$
is assumed to be continuous for $s\in[0,1]$ with the convention that $g(0)=0$ if $h(t)\to +\infty$ as $t\to +0$.
Following Sacks and Ylvisaker (1966), we define \textit{asymptotic optimality} of a sequence of sampling designs $T_n^*$ by
$$
\lim_{n\rightarrow\infty}||X- H_{q,k}(X,T_n^*)||_p\Big{/}\inf_{T\in D_n}||X- H_{q,k}(X,T)||_p=1.
$$
Now we introduce  the classes of processes used throughout the paper.  We say that:\\
(i) $X\in\ccal^{m,\beta}([a,b],V(\cdot))$ if $X\in\ccal^m([a,b])$ and $X^{(m)}$ is \textit{locally H\"{o}lder continuous}, i.e.,  if
for all $t,t+s,\in[a,b]$,
\beq \label{sls}
    \norm{X^{(m)}(t+s)-X^{(m)}(t)} \leq V(\bar{t})^{1/2} |s|^\beta,   0<\beta\le 1,
\eeq
for a positive continuous function $V(t), t\in [a,b]$, and some $\bar{t}\in [t,t+s]$.
In particular, if $V(t)=C,\,t\in[a,b]$, where $C$ is a positive constant, then $X^{(m)}$ is \textit{H\"{o}lder continuous}, and we denote it by
$X\in\ccal^{m,\beta}([a,b],C)$
\\
(ii) $X\in\bcal^{m,\beta}([a,b],c(\cdot))$ if $X\in\ccal^m([a,b])$ and there exist $0<\beta\leq 1$ and a positive continuous function $c(t),\,t\in[a,b]$, such that $X^{(m)}$ is 
\textit{locally stationary}, i.e.,
\beq \label{cls}
    \lim_{s\to 0} \frac{\norm{X^{(m)}(t+s)-X^{(m)}(t)}}{|s|^\beta}=c(t)^{1/2}
    \mbox{ uniformly in } t\in [a,b],
\eeq
(iii)
$X\in \ccal\bcal^{m,\beta}((0,1],c(\cdot),V(\cdot))$
if there exist $0<\beta\leq 1$ and positive continuous functions $c(t),V(t),\, t\in (0,1]$, such that $X\in \ccal^{m,\beta}([a,b],V(\cdot)) \cap \bcal^{m,\beta}([a,b],c(\cdot))$ for any $[a,b]~\subset~(0,1]$. \smallskip\\ 
By definition, we have that $V(t)\ge c(t), t\in (0,1]$.
 Moreover, if $X \in\ccal^{(m+1)}((0,1])$, then $c(t)=\norm{X^{(m+1)}(t)}^2$
 and we may set also $V(t)=c(t)$, $t\in(0,1]$. Some properties of locally stationary processes are considered in Berman (1974), H\"{u}sler (1995), and Seleznjev (2000).
\smallskip

\noindent \textbf{Example 1} Let $B(t),t\in[0,1]$, be a fractional Brownian motion with the covariance
function
$
  r(t,s)= (|t|^{2\bt}+|s|^{2\bt}-|t-s|^{2\bt})/2,
$
where $0<\bt\le 1$ is a Hurst parameter. Define a time changed version of the process~$B$,
$$
X(t):=B(\sqrt{t}),t\in[0,1].
$$
Then
$$
||X(t+s)-X(t)||=|\sqrt{t+s}-\sqrt{t}|^{\beta}\leq|s|^{\beta/2} \quad \mbox{for all }t,t+s\in[0,1].
$$
However on the nonsingular part of the domain, we get
$$
||X(t+s)-X(t)||= (\sqrt{t+s}+\sqrt{t})^{-\beta} |s|^{\beta} \qquad \mbox{for all }t,t+s\in(0,1]
$$
and therefore
$
X\in \ccal^{0,\beta/2}([0,1],1) \cap \ccal\bcal^{0,\beta}((0,1],c(\cdot),V(\cdot))
$
with $V(t)=c(t)=(4t)^{-\beta}$.
 \smallskip

For the
  differentiable  case, a representation of the approximation error
 can be obtained  as a consequence  of the following proposition (see Seleznjev, 2000),
  which is a q.m.\
  variant of  the well-known  Peano kernel theorem (see, e.g., Davis, 1975, p.\ 69).
   We formulate this proposition for further references for a closed interval $[a,b]\subseteq [0,1]$.
  Define linear operators of the following type over
  $\ccal^{m}([a,b])$,
  \beq \label{funcl}
   R(X):=\sum_{i\le m}\Bigl( \int_a^b X^{(i)}(s)d_i(s)\;ds +
   \sum_{j\le k_i}b_{ij}X^{(i)}(s_{ij})\Bigr),
  \eeq
 where the integrals and derivatives are taken in quadratic mean. The functions $d_i(\cdot)$ are
  assumed to be piecewise continuous over $[a,b]$ and all  points
  $s_{ij}\in [a,b]$. Further,  for any $s\in [a,b] $ and a
  given $m$,
  let $(t)_+:=\max(0,t)$, and
  $
    p_{s,m}(t):=(t -s)^m_+,\; t\in [a,b]
    $, be the function of $t$ for a given $s$.
  \begin{proposition}\label{pkern}{\em (q.m.\  Peano kernel theorem)\ }
  Let $R$ be of the form (\ref{funcl}) and $R(Y)=0$ if $Y\in \pcal_m([a,b])$. Then for all
  $X\in \ccal^{m+1}([a,b])$,
  \beq \label{epkern}
   R(X)= \int_a^b X^{(m+1)}(s) K_m(s)\; ds,
   \eeq
   where  {\em the Peano kernel}
   $
    K_m(s):= R(p_{s,m})/m!$, $ s\in [a,b].
    $
     \end{proposition}

     \noindent \textbf{Example 2}   For a given  point $t$ in  $ [0,1]$, denote by
 $K_{m,k}(t,s)$  the Peano kernel (as a function of $s$) of Hermite interpolation of $X\in \ccal^{m+1}([a,b])$ by $H_{k}(t)$.
 In particular, for  the {\it  two-point}
  piecewise linear interpolation of $ X \in
 \ccal^1([0,1])$ on the interval $[0,1]$,
 we have $K_{0,1}(t,s)=(t -s)_+ -t,\; t,s \in [0,1]$.\smallskip

\noindent For a zero mean fractional Brownian  motion $B(t)$,
 $t\in[0,1]$, with Hurst pa\-ra\-me\-ter $\beta$
 the $m$-fold integrated Brownian  motion $B_m(t)$ is defined  by setting
\beq \label{eq:fBm}
   B_m(t):=\frac{1}{(m-1)!} \int_0^1 (t-s)^{m-1}_+ B(s) \;ds,
    \;\; m\ge 1,\,
B_0(t)=B(t),
\eeq
$B^{(j)}_m(0)=0$, $j= 0,\ldots, m-1$, $B^{(m)}_m(t)=B(t).$
 Let $R_{k,2}(B_m)(t)$  be the remainder for the
 two-point Hermite
interpolation of $B_m(t),\,t\in[0,1]$, by the Hermite spline
$H_{k}(B_m)(t)$ with the norm
$$
b^{m,\beta}_{k,p}:=\norm{R_{k,2}(B_m)}_p.
$$
Explicit expressions for particular values of $b^{m,\beta}_{k,p}$ can be found in Seleznjev (2000).\\

Recall that a positive function $f$ is called \emph{regularly varying }(on the right) at the origin with index $\rho$ if for all $\lambda>0$,
\beq \label{eq:RV}
\frac{f(\lambda x)}{f(x)}\rightarrow \lambda^\rho \quad  \mbox{ as } x \to 0+.
\eeq
We denote this property by $f\in \rcal_{\rho}(0+)$. Throughout the proofs we use the following properties of
regularly varying functions  (see, e.g., Bingham et al., 1987):\\
%\citep[see, e.g.,][]{Bingham1987}:
\textnormal{(R1)} convergence in (\ref{eq:RV}) is uniform for all intervals $0<a\le \lambda \le b<\infty$;\\
% \item[\textnormal{(R2)}] $f(x)\in \rcal_{\rho}(0+)$ iff  $f(x)=x^{\rho}u(x)$, where $u(x)$ is  \emph{a slowly varying function}, i.e.,
%  $u\in \rcal_{0}(0+)$;
\textnormal{(R2)} if $f(x)\in \rcal_{\rho}(0+)$ and  $F(x)=\int_{0}^{x}f(v)dv,\;  x\in[0,1]$, then  $F\in \rcal_{\rho+1}(0+).$\\
Write $f \in \rcal^+(r(\cdot),0+)$, if for some $r\in \rcal_{\rho}(0+)$, $a>0$, we have $f(t)\le r(t)$ for all $t\in (0,a]$.

\section{Results}\label{se:Res}

We consider first the case when an optimal approximation rate  can be attained by using a suitable composite splines and sequences of sampling
designs (or generating densities).
Let $X~\!\!\!\in\ccal^{l,\alpha}([0,1],M)\cap\ccal\bcal^{m,\beta}((0,1],c(\cdot),V(\cdot))$, $l+\alpha\leq m+\beta$. We formulate the following condition for a local H\"{o}lder function $V$
and a sequence generating density $h$:\\
(C) let $g \in \rcal^+(r(\cdot),0+)$, where
\beq \label{th:RvMain}
    r(s)=\emph{{o}}(s^{(m+\beta)/(l+\alpha+1/p)-1}) \szer;
\eeq
 if $p=\infty$, then $V(t)^{1/2}r(H(t))^{m+\beta}\to 0$ as $t\to 0$;\\
 if $1\leq p<\infty$  and, additionally,  $V(G(\cdot))^{1/2}\in \rcal^+(R(\cdot),0+)$, then
 $R(H(t))r(H(t))^{m+\beta}~\in~L_p[0,b]$  for some $b>0$.
\smallskip

\noindent In the following theorem, we describe the class of generating densities eliminating the effect of the singularity point for the asymptotic
approximation accuracy.
\begin{theorem} \label{th:Main}
Let $X\in  \ccal^{l,\ah}([0,1],M)\cap \ccal\bcal^{m,\beta}((0,1],c(\cdot),V(\cdot))$, $l+\alpha\leq m+\beta$,
with the mean
  $f\in C^{m,\theta}([0,1],C)$, $\bt<\theta\le
  1$,
 be interpolated by a composite Hermite spline ${H}_{q,k}(X,T_n)$,
 $l\leq q$, $m\leq k$, where $T_n$ is a qRS($h$).
Let for the density $h$ and the local H\"{o}lder function $V$, the condition \textnormal{(C)} hold.
% and also
%\beq \label{th:RvMain}
%    r(s)=\emph{{o}}(s^{(m+\beta)/(l+\alpha+1/p)-1}) \szer.
%\eeq
Then
\beq \label{th:As}
\lim_{n\rightarrow\infty}n^{m+\beta}||X-{H}_{q,k}(X,T_n)||_p = b^{m,\beta}_{k,p} ||c^{1/2}h^{-(m+\beta)}||_p >0.
\eeq
\end{theorem}
 \noindent \textbf{Remark 1 }   (i) The condition (C) is technical and we conjecture that the main result of the above theorem is valid for a wider class
 of processes and generating functions. In particular, for $p=\infty$, if in (C) the following condition is used,
 \beq \label{eq:Cvar}
V(t)^{1/2}r(H(t))^{m+\beta}\to D \tzer \eeq for some constant $D\ge 0$, it follows straightforwardly from the proof that
\beq \label{eq:Order}
||X-{H}_{q,k}(X,T_n)||_{\infty}\asymp n^{-(m+\beta)} \mbox{ as } n\to\infty.
\eeq
More general conditions  are considered  also in Proposition \ref{pr:MainIntermed}.\\
(ii)
Let $\gamma:=1/(m+\beta+1/p)$  and  write
 $$h^*(t):=c(t)^{\gamma/2}\Big{/}\int_0^{1}c(s)^{\gamma/2}ds,\quad t\in (0,1].$$
It is well known that $h^*$ minimizes the asymptotic constant  in (\ref{th:As}) (see,~e.g.,~Seleznjev,~2000). Thus if the conditions of the
above theorem hold for $X$ and $h^*$, then $h^*$ is the asymptotically optimal density and
$$
\lim_{n\rightarrow\infty}n^{m+\beta}||X- H_{q,k}(X,T_n(h^*))||_{p}=b_{k,p}^{m,\beta}\;||c^{1/2}||_{\gamma}.
$$
For example, if  $c(t)\sim C t^{-{2\theta}}$ for some constant $C>0$ and $V(t) \asymp c(t)$  as  $t \to 0$,
then for any $1\le p<\infty$ the inequality
  \begin{equation}\label{thcond}
   \left(1-\frac{l+\alpha+1/p}{m+\beta}\right)\Big{/}\gamma< \theta <\frac{1}{\gamma} .
  \end{equation}
is a sufficient condition for $h^*$ to be the asymptotically optimal density, as it implies  the condition (C).
However for $p=\infty$, only (\ref{eq:Cvar}) holds and it follows from (\ref{eq:Order}) that the approximation error is of order $n^{-(m+\beta)}$ for large $n$.
\smallskip

\noindent\textbf{Example 3} Consider $X(t),t\in[0,1]$, defined in Example 1. By previous calculations we have $c(t)= (4t)^{-\beta}, 0<\beta<1$, and $\gamma=1/(\beta+1/p)$.
Now the condition \eqref{thcond}, i.e.,
$$
\left(1-\frac{{\beta}/{2}+{1}/{p}}{\beta}\right)\left(\beta+\frac{1}{p}\right)
=\frac{\beta}{2}- \left(\frac{1}{2p}+\frac{1}{\beta p^2}\right)
<
\frac{\beta}{2}
<
{\beta}+\frac{1}{p},
$$
is satisfied for any $1\leq p<\infty$, hence
$$
h^{*}(t)=\left(1-{\beta\gamma}/{2}\right)t^{-\beta\gamma/2} ,\qquad t\in(0,1],
$$
is the asymptotically optimal density.

\medskip
We proceed to some cases not included in Theorem \ref{th:Main} but important for applications. Specifically,
 like in the conventional approximation
for interpolation by piecewise polynomials (or splines) of order $k$ and smoothness $m>k$, the upper bound for the optimal  approximation rate
is $n^{-(k+1)}$.
We investigate the corresponding generating densities leading to the asymptotically optimal solution in such setting.
Namely, let $X\in \ccal^{l,\alpha}([0,1],M)\cap \ccal^{m}((0,1])$ be approximated by a composite Hermite spline $H_{q,k}$.
  For any $q \leq k< m $ we have  $c_k(t):=||X^{(k+1)}(t)||^2$.
We introduce  the following modification of  the condition (C), where $k+1$ is used instead of $m+\beta$:\\
\noindent(C$'$) let $g \in \rcal^+(r(\cdot),0+)$, where
\beq \label{th:RvMain2}
    r(s)=\emph{{o}}(s^{(k+1)/(l+\alpha+1/p)-1}) \szer;
\eeq
 if $p=\infty$, then $c_k(t)^{1/2}r(H(t))^{k+1}\to 0$ as $t\to 0$;\\
 if $1\leq p<\infty$  and, additionally,  $c_k(G(\cdot))^{1/2}\in \rcal^+(R(\cdot),0+)$, then
 $R(H(t))r(H(t))^{k+1}\in L_p[0,b]$  for some $b>0$.
\smallskip

\begin{theorem} \label{th:MainA1}
Let $X\in \ccal^{l,\ah}([0,1],M) \cap \ccal^{m}((0,1])$, $l+\alpha\leq m$,
with the mean $f\in C^{m+1,\theta}([0,1],C)$, $0<\theta\le 1$,
 be interpolated by a composite Hermite spline ${H}_{q,k}(X,T_n)$,
 $l\leq q$, $k<m$, where $T_n$ is a qRS($h$).
Let for the density $h$ and the function $V$, the condition \textnormal{(C$'$)} hold.
%  and also
%\beq \label{th:RvMain2}
%    r(s)=\emph{{o}}(s^{(k+1)/(l+\alpha+1/p)-1}) \szer.
%\eeq
Then
\beq \label{th:As2}
\lim_{n\rightarrow\infty}n^{k+1}||X-{H}_{q,k}(X,T_n)||_p = b^{k,1}_{k,p}\; ||c_k^{1/2}h^{-(k+1)}||_p >0,
\eeq
where $b_{k,p}^{k,1}=B(p(k+1)/2+1,p(k+1)/2+1)^{1/p}/(k+1)!$, $1\leq p<\infty$, and $b_{k,\infty}^{k,1}=2^{-(k+1)}/(k+1)!$, and $B(\cdot,\cdot)$ denotes the beta function.
\end{theorem}
Similarly to the above case $m\le k$, we consider an optimal density for such spline approximation.
\smallskip

\noindent \textbf{Remark 2 }\label{re:Main2}
Let $\gamma_k:=1/(k+1+1/p)$  and  write
 $$h^*(t):=c_k(t)^{\gamma_k/2}\Big{/}\int_0^{1}c_k(s)^{\gamma_k/2}ds,\qquad t\in (0,1].$$  If the conditions of Theorem \ref{th:MainA1} hold
 for $X$ and $h^*$, then $h^*$ is an asymptotically optimal density and
$$
\lim_{n\rightarrow\infty}n^{k+1}||X- H_{q,k}(X,T_n(h^*))||_{p}=b^{k,1}_{k,p}\; ||c_k^{1/2}||_{\gamma_k}.
$$
\smallskip

\noindent The following proposition says that an improvement  of the approximation rate is gradual with respect to
 some properties of knot generating densities.
 We introduce  conditions for a more general class of processes and generating
 densities when approximation rates $n^{-\kappa}$, $\kappa\in[l+\alpha,m+\beta]$, are attained.

\begin{proposition}\label{pr:MainIntermed}
Let $X\in \ccal^{l,\ah}([0,1],M)\cap \ccal^{m,\beta}((0,1],V(\cdot))$,  $l+\alpha \leq m+\beta$, with nonincreasing $V(\cdot)$,
be interpolated by a composite Hermite spline ${H}_{q,k}(X,T_n)$,
 $l\leq q,\, m\leq k$, for  sampling designs $T_n$ such that
\begin{equation}\label{eq:MainIntermed}
\begin{aligned}
 h_1&\le M^{-1/(l+\alpha+1/p)} n^{-\kappa/(l+\alpha+1/p)},\\
 h_j&\le V({t_{j-1}})^{-1/(2(m+\beta+1/p))} n^{-\kappa/(m+\beta+1/p)},\quad j=2,\ldots,n,
\end{aligned}
\end{equation}
where  $\kappa \in [l+\alpha,m+\beta]$. Then
$$
n^{\kappa} ||X-{H}_{q,k}(X,T_n)||_{p} \le C
$$
for some positive $C$.
\end{proposition}

\noindent \textbf{Remark 3}  The constant $C$ in the above proposition depends only on the corresponding Peano kernels of the composite spline and hence
some conventional (deterministic) approximation theory results can be used (see, e.g., Davis, 1975, p.\ 69).

\section{Numerical experiments} \label{se:Exp}
In this section, we present some examples illustrating the obtained results.
For given  knot densities and covariance functions,  first the pointwise approximation errors are found  analytically.
Then numerical maximization and integration are used to evaluate the approximation errors for the whole interval.
Notice that the nonlinear approximation here gives a
significant gain in the order of approximation error.
We consider the following knot densities
$$
h_\lambda(t)=\frac{1}{\lambda}t^{\frac{1}{\lambda}-1},\quad t\in (0,1], \quad \lambda>0,
$$
say, {\it power densities}. This leads to the sampling points
$
t_i=\left({i}/{n}\right)^\lambda, i=0, \ldots, n$, where   $\lambda=1$ corresponds to the uniform knot distribution (equidistant sampling).
Define by
$$
\dt_n(X,H_{q,k},h)(t):=X(t)-H_{q,k}(X,T_n(h))(t),\qquad t\in[0,1],
$$
 the deviation process for the approximation of $X$ by the composite spline $H_{q,k}$ with $n$ knots generated by the
density $h$ and write
$$
e_n(h)=e_n(h,X,H_{q,k},p):=||\dt_n(X,H_{q,k},h)||_p
$$
for the corresponding mean error with the density function $h$.
\smallskip

\begin{figure}[thb]
\begin{center}
\begin{tabular}{cc}
\includegraphics[height=1.6in]{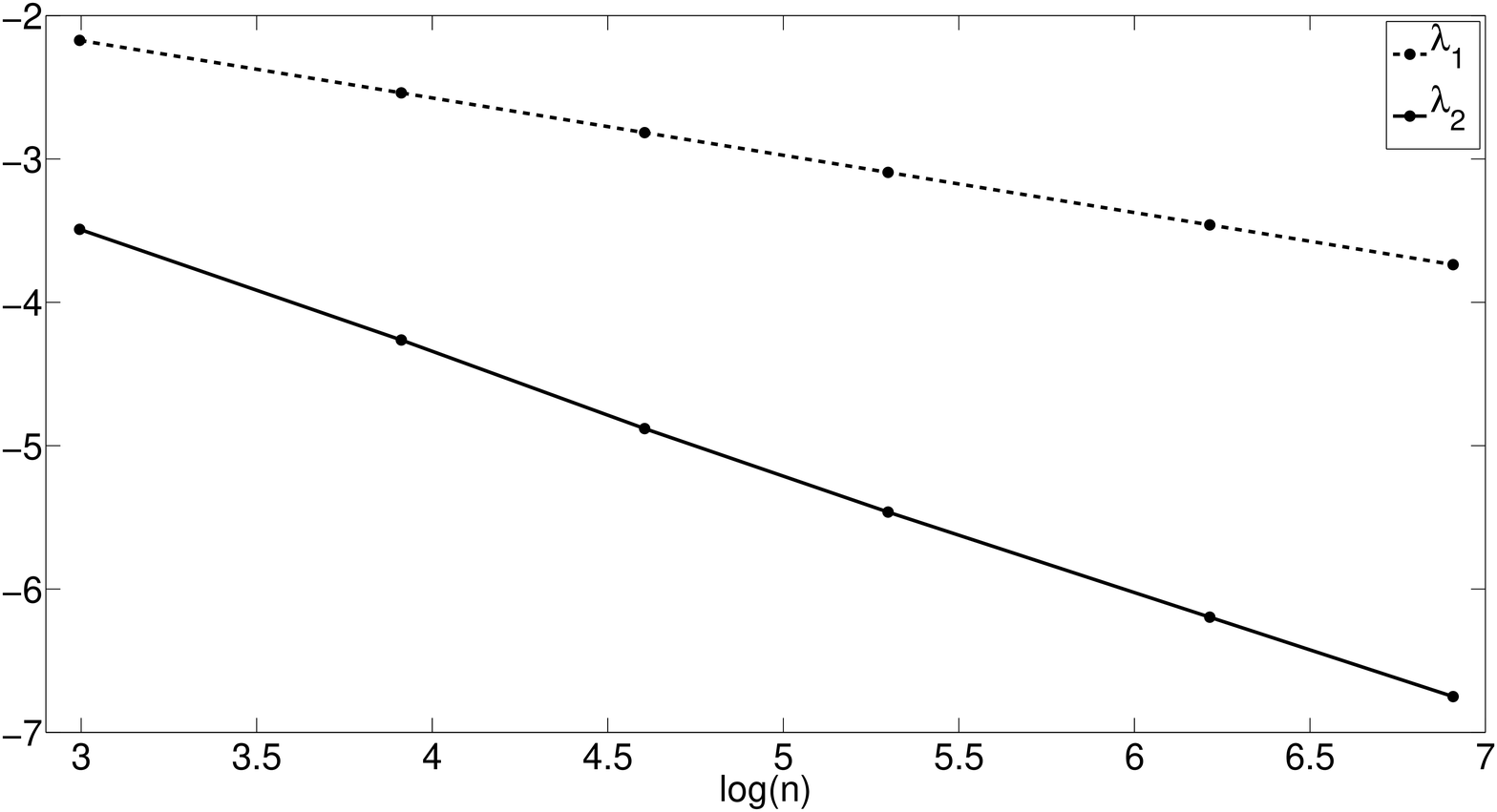} &
\includegraphics[height=1.6in]{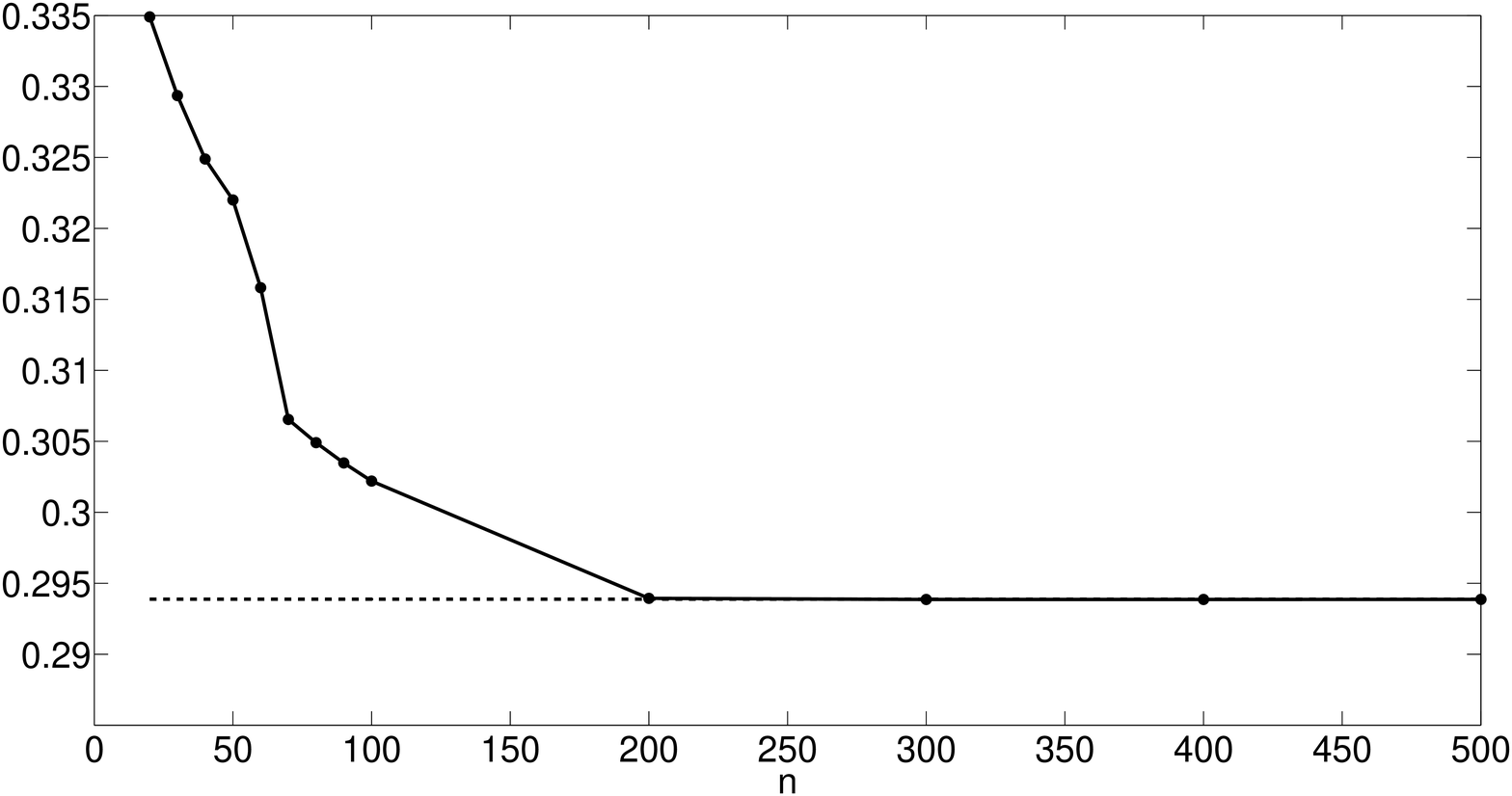}\\
\small{(a)} & \small{(b)}
\end{tabular}
\end{center}
\caption{(a) The (fitted) plots of $e_n(h_{\lambda_{1}})$ (dashed-line) and $e_n(h_{\lambda_{2}})$ (solid line) versus $n$ in a log-log scale.
(b) The convergence of $n^{0.8}e_n(h_{\lambda_{2}})$ (solid line) to the asymptotic constant (dashed line).}
\label{fg:Mx}
\end{figure}

\noindent\textbf{Example 4}
Consider $X(t),t\in[0,1]$, defined in Example 1 and let $\beta=0.4$. Let $e_n(h_{\lambda})=e_n(h_{\lambda},X,H_{1,1},\infty)$ be the mean maximal approximation error,
 $p=\infty$, when the both splines are piecewise linear interpolators, $q=k=1$.
For the power density $h_\lambda$ and $X$,
the conditions of Theorem~\ref{th:Main}\- are sa\-ti\-sfied if $\lambda>2$.
We choose the following values of the power density parameter: $\lambda_{1} = 1$ (uniform knot distribution) and $\lambda_{2} = 2.1$.
Figure \ref{fg:Mx}(a) shows the (fitted) plots and evaluated values of the mean maximal errors  $e_n(h_{\lambda_{i}}), i=1,2$, versus $n$ (in a log-log scale).
 These plots correspond to the following asymptotic behavior of the approximation errors:
$$
\begin{array}{rcl}
 e_n(h_{\lambda_{1}})&\sim& C_1 \, n^{-0.4},\quad   C_1 \simeq 0.377,\\
 e_n(h_{\lambda_{2}})&\sim& C_2 \, n^{-0.8}, \quad   C_2 \simeq  0.295 \ninf.
\end{array}
$$
For example,  the minimal number of observations needed to obtain the accuracy $0.01$ is approximately $8727$ for the equidistant sampling density $h_{\lambda_1}$, whereas it needs  only
$69$ knots when $h_{\lambda_2}$ is used, i.e., Theorem \ref{th:Main} is applicable.
Figure \ref{fg:Mx}(b) demonstrates the convergence of the scaled approximation error $n^{0.8}e_n(h_{\lambda_{2}})$ to the asymptotic constant obtained in Theorem~\ref{th:Main}.

\begin{figure}[htb]
\begin{center}
\begin{tabular}{cc}
\includegraphics[height=1.6in]{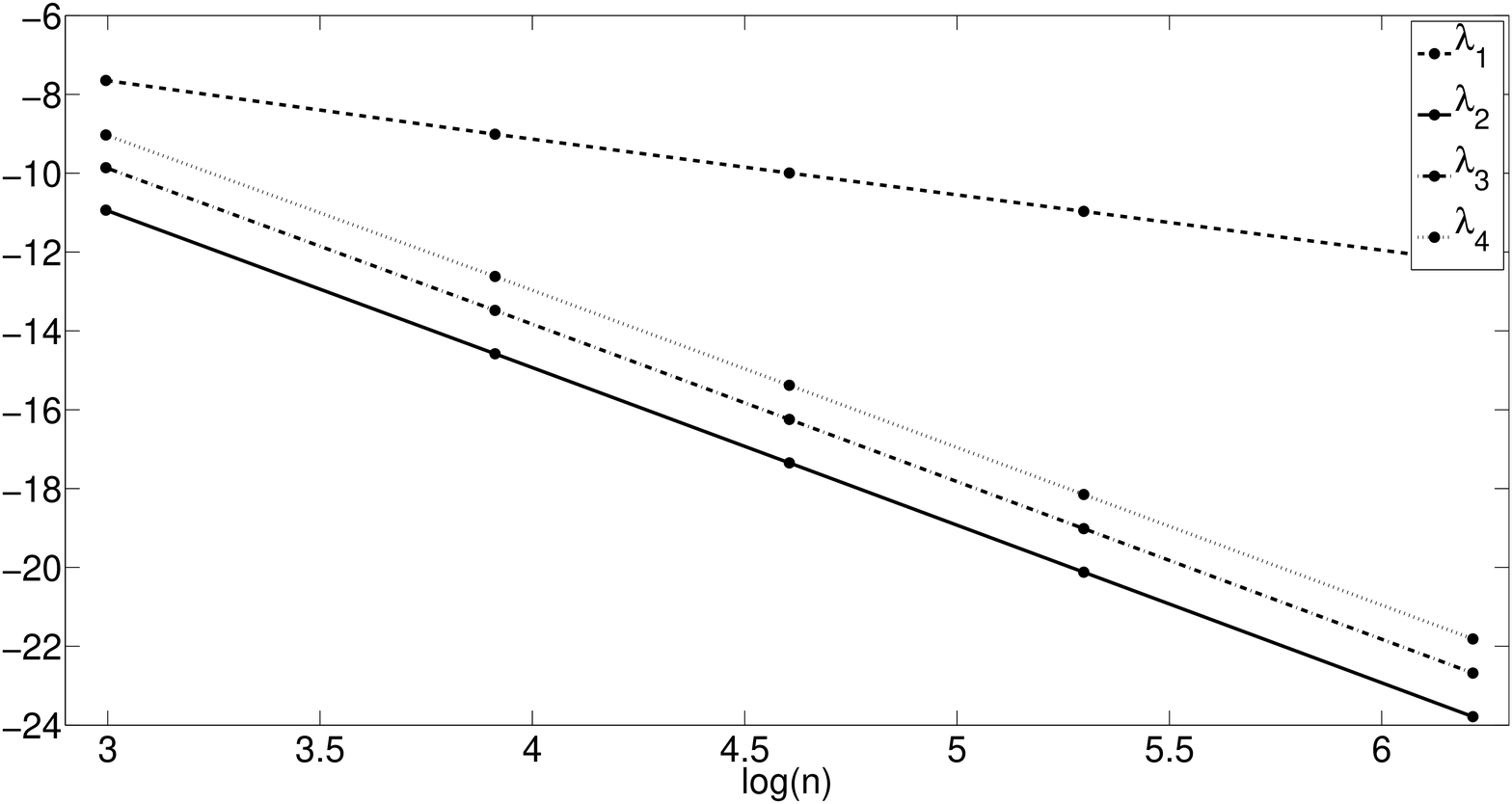} & \includegraphics[height=1.6in]{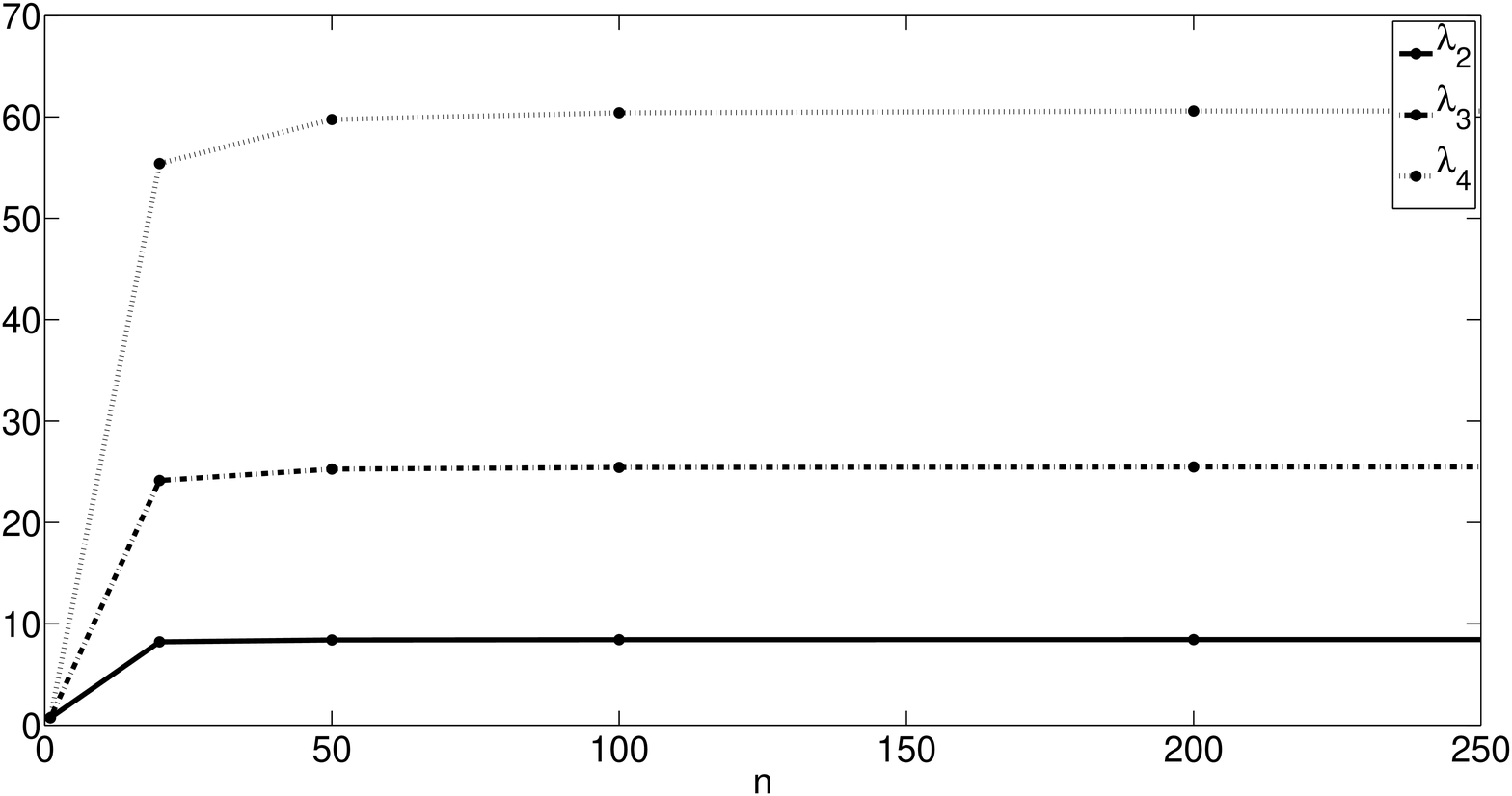}\\
\small{(a)} & \small{(b)}
\end{tabular}
\end{center}
\caption{
(a) The (fitted) plots of  $e_n(h_{\lambda_{i}}),\, i=1,2,3,4$ versus $n$ in a log-log scale.
(b) The (fitted) plots of ratios  $e_n(h_{\lambda_{i}})/e_n(h^*),\, i=2,3,4$, for  the asymptotically optimal density $h^*$.
}
\label{fg:Sm}
\end{figure}

\medskip
\noindent\textbf{Example 5}
Let $Y(t),\,t\in[0,1]$, be a zero mean stationary Gaussian process with the covariance function
$
r(s,t)=\exp\{-(s-t)^2\}.
$
We consider a distorted version of the process,
$$X(t):=t^{0.9}Y(t),\,t\in[0,1].$$
Then
$X\in \ccal^{0,0.9}([0,1],1) \cap \ccal^{m}((0,1],c_m(\cdot))$  for any $m\geq0$,
where $c_m(t)=||X^{(m+1)}(t)||^2,$ $t\in(0,1]$.
The process $Y$ has infinitely many q.m.\ derivatives, hence the approximation rate of Hermite spline approximation is limited by the order of spline only.
However, the linear methods
applied to $X$ would suffer a substantial loss of efficiency.
Consider now an approximation of $X$ by the composite Hermite spline $H_{1,3}$, i.e., a linear function on the interval adjacent to singularity, and the cubic Hermite spline otherwise.
We investigate the mean integrated ($p=2$) approximation error $e_n(h_{\lambda})=e_n(h_{\lambda},X,H_{1,3},2)$.
For the corresponding local stationarity function, we have
$$
%c(t)=1680 t^\frac{9}{5}-\frac{44631}{1250}t^{-\frac{11}{5}}-\frac{16929}{125000} t^{-\frac{21}{5}}+\frac{8424}{5} t^{-\frac{1}{5}}+\frac{4322241}{10^{8}}t^{-\frac{31}{5}}
c_3(t)\sim C t^{-6.2} \tzer,
$$
where $C\simeq 0.0432$.
For the power density $h_\lambda$ and $X$, the conditions of Theorem \ref{th:MainA1} are satisfied if $\lambda> 20/7$.
We choose $\lambda_{2}=3$,  $\lambda_{3} = 4$, and $\lambda_{4} = 5$ satisfying this condition, together with
$\lambda_{1} = 1$, the equidistant sampling density.
Figure \ref{fg:Sm}(a) shows the (fitted) plots and evaluated values of mean integrated errors  $e_n(h_{\lambda_{i}})$, $i=1,2, 3,4$, versus $n$
(in a log-log scale).

\noindent Note that the estimated approximation  rate for the uniform knot distribution is $n^{-1.4}$,  whereas   the power densities with parameters
$\lambda_2,\lambda_3$, and $\lambda_4$ result in the regular (nonsingular) case rate, namely, $n^{-4}$.
By Remark~2, the asymptotically optimal density
$$
 h^*(t)=c_3(t)^{\frac{1}{9}}\Big{/}\int_0^{1}c_3(s)^{\frac{1}{9}}ds,\quad  t\in (0,1].
$$
Moreover, $h^*(t)\sim h_{\lambda}(t)$ with $\lambda=45/14 \simeq 3.21$ as $t\to0$. Figure \ref{fg:Sm}(b) demonstrates the
ratios $e_n(h_{\lambda_{i}})/e_n(h^*),\, i=2,3,4$  and an essential benefit in the asymptotic constant for the optimal density $h^*$.

\section{Proofs} \label{se:Pro}

\noindent {\it Proof of Theorem \ref{th:Main}.}
We investigate the asy\-mpto\-tic be\-havi\-or of the q.m.\ error
$e_{n}(t) : = \norm{X(t)- H_{q,k}(X,T_n)(t)}$ for any  $t\in [t_{j-1},t_j]$, $i=1,\ldots,n$, when the number of knots $n$ tends to infinity.
An explicit formula when $l=0$ or $m=0$ is used, and the q.m.\  Peano kernel representation is applied otherwise.
We show that the densities satisfying the conditions of the theorem lead asymptotically to
 eliminating the effect of the singularity point.
Further, we find the asymptotic form  of $\norm{X- H_{q,k}(X,T_n)}_p$ for any density $h(\cdot)$ satisfying the conditions.
\smallskip

 Throughout the proof we utilize the  next property of the  Peano kernel (see, e.g.,  Seleznjev 2000):\\
 Let $K_{m-1,k}(t,s), t,s\in [0,1]$,
 be a Peano kernel for
the two-point Hermite interpolation, $m\geq 1$. If  $m\le k \le 2m+1$, then
\begin{equation} \label{pkpros}
\int_0^1 K^*_{m-1, k}(t,s) ds=0.
\end{equation}

\noindent Additionally, the following (\emph{shifting}) property of a regular varying function $r$
 is used:\\
 there exist  positive $C, \beta$ such that
\beq \label{eq:RVpr}
 r(s)\le C r(v)  \mbox{ for all  } \frac{1}{2}\le \frac{s}{v}\le {2},\qquad v\neq0,\,\,\max(s,v)<\beta,
\eeq
which follows directly from  the regular variation property (R1).

\noindent Note  that \eqref{th:RvMain} and the regular variation property (R2) imply
\begin{equation} \label{eq:Gas}
G(s)=o\left(s^{(m+\beta)/(l+\alpha+1/p)}\right)\; \mbox{ as } s\to 0.
\end{equation}

Now we come to the detailed presentation of the proof. We investigate the asy\-mpto\-tic be\-havi\-our of the q.m.\  error  $e_{n}(t) : =
\norm{X(t)-H_{q,k}(X,T_n)(t)}$ for any  $t\in [t_{j-1},t_j]$, $i=1,\ldots,n$, when the number of knots $n$ tends to infinity.
Consider first the case $p=\infty$, the mean maximal norm.
Let for a fixed $\delta>0$, ${J_\delta}$ denote the largest integer such that $t_{J_\delta}<\delta$.
Then the error can be decomposed as follows,
$$
\begin{aligned}
%e_n:&=
||X- H_{q,k}(X,T_n)||_\infty=\max_{[0,1]}e_n(t)=\max( e_{n,j},j=1,\ldots,n) =\max(S_1,S_2,S_3),
\end{aligned}
$$
where $e_{n,j}:=\max_{[t_{j-1},t_j]}e_n(t)$,\\
$$
\begin{aligned}
S_1&:=S_1(n)=e_{n,1},\\
S_2&:=S_2(n)=\max( e_{n,j},j=2,\ldots,J_\dt),\\
S_3&:=S_3(n)=\max( e_{n,j},j=J_\dt+1,\ldots,n).
\end{aligned}
$$

For  the first interval and $S_1$, if $l=0$, then
 a piecewise linear interpolator $H_1(X,T_n)$ is applied. Hence, by the definition, (R2), (C), and (\ref{eq:Gas}), we obtain that
\beq \label{eq:In1}
\begin{aligned}
n^{m+\beta} e_{n,1}&\le M 2^{-\alpha} n^{m+\bt}2 h_1^{\alpha}\le M 2^{-\alpha}n^{m+\bt} G(1/n)^{\alpha}
\\&=n^{m+\bt}\mbox{o}(n^{-(m+\bt)})=\mbox{o}(1) \ninf.
\end{aligned}
\eeq
If $l>0$, Proposition \ref{pkern} and \eqref{pkpros}  imply, that for $t\in[0,t_1]$,
$$
X(t)-H_q(X,T_n)(t)=
h_1^l\int_0^1 (Z(\os)-Z(0)){K}^*_{l-1,q}(\ot,\os) \;d\os,
$$
where   $Z(\os)= \xder{l}(t_{j-1} +h_j\os) $, and ${K}^*_{l-1,q}(\ot,\os)$ denotes
the corresponding Peano kernel for the two-point Hermite interpolation at $[0,1]$. Now the H\"older property  and (\ref{eq:Gas}) give
\beq \label{eq:In11}
n^{m+\bt} e_{1,n}\leq C n^{m+\bt} h_1^{l+\alpha}= C n^{m+\bt} G(1/n)^{l+\alpha}=\mbox{ o }(1) \mbox{ as } n\to \infty.
\eeq

For $S_2$, let first $m\geq 1$. Then  Proposition \ref{pkern} implies
\beq \label{eq:Pkr1}
\begin{aligned}
X(t)- H_{q,k}(X,T_n)(t)&=X(t)- H_{k}(X,T_n)(t)
=
\int_{t_{j-1}}^{t_j} X^{(m)}(\os){K}_{m-1,k}(t,s) \;ds
\\&=h_j^{m}\int_0^1 Y(\os){K}^*_{m-1,k}(\ot,\os) \;d\os
\end{aligned}
\eeq
 where   $Y(\os)= \xder{m}(t_{j-1} +h_j\os) $, and ${K}^*_{m-1,k}(\ot,\os)$ denotes
the corresponding Peano kernel for the two-point Hermite interpolation at $[0,1]$. Using \eqref{pkpros}, we get
\beq \label{eq:Pkr2}
X(t)- H_{q,k}(X,T_n)(t)=h_j^{m}\int_0^1 \left(Y(\os)-Y(0)\right){K}^*_{m-1,k}(\ot,\os) \;d\os.
\eeq
Now the Cauchy-Schwartz inequality together with (\ref{sls}) give
\begin{equation}\label{eq:In111}
e_{n,j}\leq C h_j^{m+\bt} V(w_j)^{1/2},
\end{equation}
where
$V(w_j)^{1/2}=\max_{[t_{j-1},t_j]}V(t)^{1/2}, w_j=G(s_j)$, $s_j\in \left[{(j-1)}/{n},{j}/{n} \right]$, and
 $$
 C=C(K_{m-1,k})=\max_{t\in[0,1]}||K_{m-1,k}(t,\cdot)||_2.
 $$
Notice that by definition and the integral mean value theorem  for a positive $C_1$,
\beq \label{eq:Hj1}
h_j= G\left(\frac{j}{n}\right)-G\left(\frac{j-1}{n}\right)= g(v_j)\frac{1}{n}\le C_1\frac{1}{n}, j=1, \ldots,n,
\eeq
and, therefore, from (\ref{eq:Pkr2}) we obtain, $v_j\in \left[{(j-1)}/{n},{j}/{n} \right]$,
$$
\begin{aligned}
n^{m+\bt}e_{n,j}&\le  C  n^{m+\bt} V(w_j)^{1/2} h_j^{m+\bt}  \le  C V(w_j)^{1/2} g(v_j)^{m+\bt}\\&\le  C V(w_j)^{1/2} r(v_j)^{m+\bt}.
\end{aligned}
$$
Now the (shifting) property (\ref{eq:RVpr}) implies for some positive $C_2$,
\beqr \label{eq:S2}
 n^{m+\bt}e_{n,j}&\le&  C\, V(G(s_j))^{1/2}r(v_j)^{m+\bt} \nonumber \\
&\le& C_2 V(G(s_j))^{1/2}r(s_j)^{m+\bt} ,\; j=2,\ldots, J_\dt,
\eeqr
for  small enough $\dt$  by the definition, since $s_j, v_j \in [(j-1)/{n},{j}/{n}]$ and
$$
\frac{1}{2} \le \frac{j-1}{j}\le \frac{s_j}{v_j} \le \frac{j}{j-1}\le 2, \quad j\ge 2.
$$

For $m=0$, we have $k=1$, and  the explicit expression for the piecewise linear interpolator yields
$$
\begin{aligned}
e_{n,j}(t)&=\,\, ||X(t)-(1-\bar t) X(t_{j-1})-\bar t X(t_j)||\leq ||X(t)-X(t_{j-1})|| + \bar t ||X(t_j)-X(t_j-1)||\\
&\le (\bar t^{\bt}+ \bar t)V(w_j)^{1/2}h_j^{\bt}\le 2V(w_j)^{1/2}h_j^{\bt},
\end{aligned}
$$
and therefore, by (\ref{eq:Hj1}), the H\"older condition, and the shifting property (\ref{eq:RVpr}) we get
\begin{equation}\label{eq:S20}
\begin{aligned}
n^{\bt}e_{n,j} \leq 2 g(v_j)^{\bt}V(G(s_j))^{1/2}\leq C  r(s_j)^{\bt}V(G(s_j))^{1/2},
\end{aligned}
\end{equation}
for some positive constant $C$. Thus, by  (C), (\ref{eq:S2}), and (\ref{eq:S20}) for any $\eps>0$, sufficiently small $\dt$, and $n$ large enough, we get
\beq  \label{eq:In2}
n^{m+\beta}S_2=\max( n^{m+\bt}e_{n,j},j=2,\ldots, J_\dt) <\eps,
\eeq
for any $m\geq 0$.

For $S_3$,  we show first that
\beq
n^{m+\beta}S_{3} \to b^{m,\bt}_{k,\infty} \max_{[\dt,1]}(c(t)^{1/2}h^{-(m+\bt)}(t)) \ninf.
\eeq
The main steps of the proof repeat those of the corresponding result for a regular sequence of designs in Seleznjev (2000)
%\citet{Seleznjev2000}
 for interpolation of a smooth random function by the
Hermite spline $H_k(X,T_n)$. Applying the local stationarity condition  and uniform continuity of  positive $c(t), h(t), t\in [\dt, 1]$, we obtain
\begin{equation} %\label{eq:In3}
\begin{aligned}
n^{m+\bt} S_3=&n^{m+\bt}\max(e_n(t), t\in [t_{j-1},t_j], j=J_\dt+1, \ldots,n) \\
=& b^{m,\bt}_{k,\infty}  \max(c(w_j)^{1/2}h(w_j)^{-(m+\bt)},  j=J_\dt+1, \ldots,n)(1+\oo(1)),\\
\end{aligned}
\end{equation}
for some  $w_j\in[t_{j-1},t_j], j=J_\dt+1, \ldots,n$, and consequently
\begin{equation} \label{eq:In3}
\begin{aligned}
n^{m+\bt} S_3=& b^{m,\bt}_{k,\infty}  \max(c(t)^{1/2}h(t)^{-(m+\bt)}, t\in [t_{J_\dt},1])(1+\oo(1)) \ninf.
\end{aligned}
\end{equation}
Moreover,
\beq \label{eq:In31}
\begin{aligned}
  \max(c(t)^{1/2}h&(t)^{-(m+\bt)}, t\in [t_{J_\dt},1])\\
 &=\max(c(t)^{1/2}h(t)^{-(m+\bt)}, t\in [\dt,1])(1+\oo(1))  \ninf,
\end{aligned}
\eeq
since (\ref{eq:Hj1}) yields $|t_{J_\dt}-\dt|\le h_{J_\dt}\le C /n$ for a positive constant $C$.\\
Let $a_\dt:=\max(c(t)^{1/2} h(t)^{-(m+\bt)}, t\in [t_{J_\dt},1])$. Note that it follows from (C) that $a_\delta$ is bounded, and the
 monotone convergence gives
\beq \label{eq:Adt}
 a_\dt \uparrow a:=b^{m,\bt}_{k,\infty}  \sup(c(t)^{1/2}h(t)^{-(m+\bt)}, t\in [0,1]) \mbox{ as } \dt\to 0.
\eeq

So, first, we select  $\dt$ sufficiently small and  apply  (\ref{eq:In2}), (\ref{eq:In3}), and (\ref{eq:Adt})
for sufficiently large $n$. Then for the selected $\dt$,
(\ref{eq:In1}) and (\ref{eq:In11}) imply the assertion. This completes the proof for~$p=\infty$.\\\medskip

The proof for $1\leq p<\infty$ is analogous to the previous case, so we give the main steps only.
Let for a fixed $\delta>0$,
$$
e_n^p:=\int_{0}^{1} e_n(t)^p dt=
\sum_{j=1}^n \int_{t_{j-1}}^{t_j} e_n(t)^p dt= e_{n,1}^p+\sum_{j=2}^{J_\dt} e_{n,j}^p+ \sum_{j=J_\dt+1}^n e_{n,j}^p=S_1+S_2+S_3,
$$
where $S_1=S_1(n):=e_{n,1}$, the sum $S_3=S_3(n)$ includes all terms $e_{n,j}$ such that $[t_{j-1},t_j]\subset [\dt, 1]$,
say, $j\ge J_\dt+1$, and $S_2=e_n^p-S_1-S_3$.

For $S_1$, the first interval, if $l=0$,  $H_1(X,T_n)$ is applied, hence by the definition, (R2), (C), and (\ref{eq:Gas}), we have
\beq \label{eq:In1p}
\begin{aligned}
n^{p(m+\bt)} (e_{n,1}(t))^{p}\le& n^{p(m+1)} M 2^{-p\alpha}\int_{0}^{h_1} 2 h_1^{p\alpha} dt\le M 2^{1-p\alpha}n^{p(m+\bt)} G(1/n)^{p\alpha+1}\\ =
&n^{p(m+\bt)}\mbox{o}(n^{-p(m+\bt)})=\mbox{o}(1) \ninf.
\end{aligned}
\eeq
If $l>0$,  Proposition \ref{pkern} and \eqref{pkpros} imply that for $t\in[0,t_1]$,
$$
X(t)-H_q(X,T_n)(t)=
h_1^l\int_0^1 (Z(\os)-Z(0)){K}^*_{l-1,q}(\ot,\os) \;d\os
$$
where   $Z(\os)= \xder{l}(t_{j-1} +h_j\os) $, and ${K}^*_{l-1,q}(\ot,\os)$ denotes
the corresponding Peano kernel for the two-point Hermite interpolation at $[0,1]$. Now the H\"older inequality  and (\ref{eq:Gas}) give
\beq\label{eq:In11p}
n^{p(m+\bt)} (e_{n,1})^{p}\leq C n^{p(m+\bt)} h_1^{1+pl+p\alpha}= C n^{p(m+\bt)} G(1/n)^{1+pl+p\alpha}=\mbox{o}(1).
\eeq

For $S_2$ and $m\geq 0$, i.e., $j=2,\ldots, J_\dt$, (\ref{eq:Pkr2}) and  (\ref{eq:Hj1}) yield that for some positive constant $C$
$$
 n^{(m+\bt)}e_{n,j}\le C n^{(m+\bt)} V(w_j)^{1/2}h_j^{m+\bt+1/p} \le  C V(w_j)^{1/2}g(v_j)^{(m+\bt)}h_j^{1/p},
$$
where $w_j,v_j\in [{(j-1)}/{n},{j}/{n} ]$.
The (shifting) property (\ref{eq:RVpr}) together with condition (C) imply that for some positive $C_1$,
\beqr \label{eq:In2p}
 n^{(m+\bt)p}e_{n,j}^p(t)&\le&  C^p V(w_j)^{p/2}g(v_j)^{(m+\bt)p}h_j \le  {C_1} R(s_j)^p r(s_j)^{(m+\bt)p}h_j\nonumber\\
 &\le&  C_1  \int_{t_{j-1}}^{t_{j}} R(H(t))^p r(H(t))^{(m+\bt)p} dt,\; j=2,\ldots, J_\dt,
\eeqr
where for $s_j=H(u_j)\in \left[(j-1)/n,{j}/{n} \right]$,
 $$
\begin{aligned}
 R(s_j)r(s_j)^{m+\bt}= R(H(u_j))r(H(u_j))^{m+\bt}=\min_{[t_{j-1},t_j]} R(H(t))r(H(t))^{m+\bt}.
\end{aligned}
$$
Thus, for any $\eps_1>0$, sufficiently small $\dt$, and large enough $n$, we obtain
\beq\label{S2p}
n^{(m+\bt)p} S_2^p\le C_1 \int_{t_{1}}^{t_{J_\dt}} R(H(t))^p r(H(t))^{(m+\bt)p} dt  <\eps_1
\eeq
Analogously to the case $p=\infty$, the explicit form of piecewise linear interpolator is used to obtain (\ref{S2p}) when $m=0$.

For $S_3$, we show that
%(m+1)
%\left(b_{k,p}^{m,1}\right)

\beq
n^{(m+\bt)p}S_{3}^p\to  \left(b_{k,p}^{m,\bt}\right)^p\int_{\dt}^1 c(t)^{p/2}h(t)^{-(m+\bt)p}dt \ninf.
\eeq
Applying the local stationarity condition and uniform continuity of positive $c(t), h(t), t\in  [\dt, 1]$, we get
$$
n^{(m+\bt)p}S_{3}^p=    \left(b_{k,p}^{m,\bt}\right)^p \left(\sum_{j=J_\dt+1}^n c(w_j)^{p/2} h(w_j)^{{-(m+\bt)p}}h_j\right)
\left(1+\oo(1)\right)
$$
 for some  $w_j\in[t_{j-1},t_j], j=J_\dt+1, \ldots,n$, and consequently
\beqr\label{eq:In3p}
n^{(m+\bt)p}S_{3}^p&=&  \left(b_{k,p}^{m,\bt}\right)^p\left(\int_{t_{J_\dt}}^1 c(t)^{p/2} h(t)^{{-(m+\bt)p}}dt\right)\left(1+\oo(1)\right)\nonumber\\
&=& \left(b_{k,p}^{m,\bt}\right)^p\left(\int_{\dt}^1 c(t)^{p/2}h(t)^{{-(m+\bt)p}}dt\right)(1+\oo(1))\nonumber\\
&=:&a_\dt(1+\oo(1))  \ninf.
\eeqr
Now it follows by (C), that  $c(t)^{p/2}h(t)^{{-(m+\bt)p}}$ is integrable on $[0,1]$, hence
\beq \label{eq:Adtp}
 a_\dt \uparrow a:=\left(b_{k,p}^{m,\bt}\right)^p\int_{0}^1 c(t)^{p/2}h(t)^{{-(m+\bt)p}}dt, t\in [0,1]) \mbox{ as } \dt\to 0.
\eeq
This completes the proof.
\medskip

\noindent\textit{Proof of Theorem \ref{th:MainA1}. }
The proof repeats the steps of that of Theorem \ref{th:Main} for the differentiable case. For the calculations of the asymptotic constants we refer to Seleznjev (2000).
\medskip

\noindent{\em Proof of Proposition \ref{pr:MainIntermed}.}
The arguments are similar for the first interval and the other intervals in both the mean integrated, and maximal norms. So we demonstrate the main
steps for the $[t_{j-1},t_j], j\ge 2$, the mean maximal norm, and the differentiable case.
Applying Proposition \ref{pkern} and
\eqref{pkpros} to the Hermite spline approximation, we have, $t\in
[t_{j-1},t_j], j\ge 2$,
$$
 X(t)-H_{k}(X)(t) = h_j^{m}\int_0^1 \left(Y(\os)-Y(0)\right){K}^*_{m-1,k}(\ot,\os) \;d\os,
$$
where   $Y(\os)= \xder{m}(t_{j-1} +h_j\os) $, and ${K}^*_{m-1,k}(\ot,\os)$ denotes
the corresponding Peano kernel for the two-point Hermite interpolation at $[0,1]$. Now the Cauchy-Schwartz inequality, \eqref{sls}, and the
monotonicity of $V$ give
\begin{equation*}
||X(t)-H_{k}(X)(t)||_{\infty}\leq C h_j^{m+\bt} V(t_{j-1})^{1/2},
\end{equation*}
where $C=C(K_{m-1,k})=\max_{t\in[0,1]}||K_{m-1,k}(t,\cdot)||_2$. Finally, applying the condition \eqref{eq:MainIntermed} implies the result.
\medskip

\medskip

\medskip

\noindent{\bf References}
\begin{reflist}

Abramowicz, K.\ and Seleznjev, O.\  (2008). On the error of the Monte Carlo pricing method. \emph{J.\ Num. and Appl.\ Math.} 96, 1-10

 Benhenni, K.\ and  Cambanis, S.\ (1992). Sampling designs for estimating
integrals of stochastic process. \emph{Ann. Statist.}  20, 161-194.

Berman, S.M.\    (1974). Sojourns and extremes of Gaussian process. {\it  Ann.\ Probab.}
     2, 999- 1026; corrections  8, 999 (1980);  12, 281,    (1984).

Bingham, N.H., Goldie, C.M., and Teugels, J.L.\ (1987). \emph{Regular variation.} Cambridge Univ. Press.

    de Boor, C.\ (1973). Good approximation by splines with variable
    knots, In: Meir, A.\     and Sharma, A., Eds., \emph{Spline functions and approximation theory},
     Birkh\"{a}user, 57-72.

Buslaev, A.P.\ and  Seleznjev, O.\ (1999). On certain extremal problems in theory of approximation of random processes. \emph{ East J.\   Approx.} 5, 467--481.

Cohen, A.  and  D'Ales, J.-P.\ (1997). Nonlinear approximation of random functions. \emph{SIAM J.\ Appl.\ Math.} 57 , 518-540.

 Cohen, A.,  Daubechies, I.,  Guleryuz, O.G.,   and Orchard, M.T.\  (2002). On the
importance of combining wavelet-based nonlinear approximation with
coding strategies. \emph{IEEE Trans. Inform. Theory} 48 , 1895-1921.

Creutzig, J.\ and  Lifshits, M.\ (2006). Free-knot spline approximation of fractional Brownian motion.
In:    Keller, A.,  Heinrich, S., and  Neiderriter, H., Eds., \emph{Monte Carlo and quasi Monte Carlo methods}. Springer, Berlin, 195-204.

Creutzig, J.,  M\"uller-Gronbach, T., and Ritter, K.\ (2007).  Free-knot spline approximation of stochastic processes.
 \emph{J.\ Complexity} 23,  867-889.

      Davis, P.J.\ (1975). \emph{Interpolation and approximation}. Dover,
New York.

DeVore, R.\ (1998). Nonlinear approximation. \emph{Acta Numer.} 7 , 51-150.

   Eplett, W.T.\ (1986). Approximation theory
for simulation of continuous Gaussian processes. {\em Prob.\ Theory Related Fields} 73,
159-181.

H\"{u}sler, J.\ (1995). A note on extreme values of locally stationary Gaussian processes.
 {\it J.\ Statist.\ Plann.\ Inference}  45, 203- 213.

H\"usler,  J.,  Piterbarg, V., and  Seleznjev, O. (2003).
 On convergence of the uniform norms for Gaussian processes and linear approximation problems.  \emph{Ann. Appl.\ Probab.} 13,  1615-1653.

 Kon, M.\ and   Plaskota, L.\ (2005).  Information-based nonlinear approximation: an average case setting. \emph{J.\ Complexity} 21,  211-229,

  Ritter,  K.\ (2000). \emph{Average-case analysis of numerical problems},
Springer-Verlag.

Sacks, J.\ and  Ylvisaker, D.\   (1966). Design for regression problems with correlated errors.
\emph{Ann.\ Math.\ Statist.} 37, 66-89.

Seleznjev, O.\ (1996).    Large deviations in the piecewise linear approximation of Gaussian
processes with stationary increments.
\emph{Adv. in\ Appl.\ Prob.} 28, 481-499.

Seleznjev, O.\ (2000).\ Spline approximation of random processes and design problems. \emph{J.\ Statist.\ Plann.\ Inference} 84, 249-262.

Su, Y.\ and Cambanis, S.\ (1993) \ Sampling designs for estimation of a random  process. \emph{Stoch. Proc.\ Appl.} 46, 47-89.

\end{reflist}
\end{document}